\documentclass[12pt]{article}

\usepackage[margin=1in]{geometry}  % set the margins to 1in on all sides
\usepackage{graphicx}              % to include figures
\usepackage{amsmath}               % great math stuff
\usepackage{amsfonts}              % for blackboard bold, etc
\usepackage{amsthm}                % better theorem environments
\newtheorem*{theorem}{Theorem}

\begin{document}

\title{Ergodic dynamical systems over the Cartesian power of the ring of p-adic integers}

\author{Valerii Sopin}

\maketitle

\begin{abstract}
For any 1-lipschitz ergodic map $F:\; \mathbb{Z}^{k}_{p} \mapsto
\mathbb{Z}^{k}_{p},\;k>1\in\mathbb{N},$ there are 1-lipschitz
ergodic map $G:\; \mathbb{Z}_{p} \mapsto \mathbb{Z}_{p}$ and two
bijection $H_k$, $T_{k,\;P}$ that $$G = H_{k} \circ T_{k,\;P}\circ
F\circ H^{-1}_{k} \text{ and } F = H^{-1}_{k} \circ
T_{k,\;P^{-1}}\circ G\circ H_{k}.$$
\end{abstract}

\section{Introduction}

The p-adic numbers, which appeared more than a century ago as a pure mathematical construction at the end of the
20th century were recognized as a base for adequate descriptions of
physical, biological, cognitive and information processing
phenomena.

Now the p-adic theory, and wider, ultrametric analysis and
ultrametric dynamics, is a rapidly developing area that finds
applications in various sciences (physics, biology, genetics,
cognitive sciences, information sciences, computer science,
cryptology, numerical methods and etc). An interested reader is referred to the monograph
[1] and references therein.

A pseudorandom number generator (PRNG) is an algorithm that takes a
short random string (a seed) and stretches it to a much longer
string that looks like random. ``Looks like random" means passes
prescribed statistical tests. Thus, the very concept of
``pseudorandomness" depends on what tests the output of the PRNG
must pass.

PRNG are being used widely: in cryptography, for computer
simulations, in numerical analysis (e.g., in quasi Monte Carlo
algorithms) and etc. However, a common demand is that the output of
a PRNG must be uniformly distributed: limit frequencies of
occurrences of symbols must be equal for all symbols.

1-lipschitz transformations on the Cartesian power of the ring of
p-adic integers can be used (and already are being used) to
construct both state transition functions and output functions of
various PRNGs.

Ergodic 1-lipschitz transformations on the Cartesian power of the
ring of 2-adic integers have been considered as a candidate to
replace linear feedback shift registers (LFSRs) in keystream
generators of stream ciphers, since sequences produced by such
function are proved to have a number of good cryptographic
properties, e.g., high linear and 2-adic complexity, uniform
distribution of subwords and etc, see [1-4].

V. Anashin gave criteria for measure-preservation and ergodicity
of 1-lipschitz transformations on the ring of p-adic integers [1-5].
Most recently V. Anashin et al [5] have used the Van der Put basis
to describe the ergodic 1-lipschitz functions on the ring of 2-adic
integers.

However, issue of describing the ergodic 1-lipschitz transformations on
the Cartesian power of the ring of p-adic integers has been opened
so far. In this paper we present the resulting
solution to this problem.

\section{P-adics}

Let $p$ be an arbitrary prime. The p-adic valuation is denoted by
$|*|_p$. We remind that this valuation satisfies the strong triangle
inequality: $$|x + y|_p \leq \max(|x|_p, |y|_p).$$

This is the main distinguishing property of the p-adic valuation
inducing essential departure from the real or complex analysis (and
hence essential difference of p-adic dynamical systems from real and
complex dynamical systems).

The ring of p-adic integers is denoted by the symbol
$\mathbb{Z}_{p}$. We remind that any p-adic integer (an element of
the ring $\mathbb{Z}_{p}$) can be expanded into the series
$$\sum\limits_{i=0}^\infty\alpha_i p^i, \text{ where } \alpha_i\in\{0,\;\dots,\;p-1\},
i\in\mathbb{N}.$$

$\mathbb{Z}^k_{p}, k\in\mathbb{N},$ is the Cartesian power of the
ring of p-adic integers. Metric on $k$-th Cartesian power
$\mathbb{Z}^k_{p}$ can be defined in a similar way:
$$|(a^1, \dots , a^k)-(b^1, \dots , b^k)|_p= \max\{|a^i -b^i |_p: i=
1, 2, \dots, k\}$$ for every $(a^1, \dots , a^k)$, $(b^1, \dots ,
b^k) \in \mathbb{Z}^k_{p}$.

The space $\mathbb{Z}^k_{p}$ is equipped by the natural probability
measure, namely, the Haar measure $\mu_{p, k}$ normalized so that
$\mu_{p, k}(\mathbb{Z}^k_{p}) = 1$.

A map $F$: $\mathbb{Z}^{k}_{p}$ $\mapsto$ $\mathbb{Z}^{k}_{p}$ is
measure-preserving, if $\mu_{p, k}(F^{-1}(S)) = \mu_{p, k}(S)$ for
any measurable set $S \subseteq \mathbb{Z}^{k}_{p}$.

A map $F$: $\mathbb{Z}^{k}_{p}$ $\mapsto$ $\mathbb{Z}^{k}_{p}$ is
ergodic, if $\mu_{p, k}(S) = 0 $ or $\mu_{p, k}(S) = 1$ follows for
any measurable set $S$: $\mu_{p, k}(F^{-1}(S)) = \mu_{p, k}(S)$.

$(\mathbb{Z}/p^n\mathbb{Z})^k$ is the Cartesian power of the residue
ring modulo $p^n$. Reduction modulo $p^n$ (it's denoted by the
symbol $\text{mod } p^n$), $n \in \mathbb{N}$, is an epimorphism of
$\mathbb{Z}^k_p$ onto $(\mathbb{Z}/p^n\mathbb{Z})^k$:
$$(\sum\limits_{i=0}^{\infty}\alpha^0_i
p^i,\sum\limits_{i=0}^{\infty}\alpha^1_{i} p^i, \dots,
\sum\limits_{i=0}^{\infty}\alpha^{k-1}_{i} p^i) \text{ mod } p^n=
(\sum\limits_{i=0}^{n-1}\alpha^0_i
p^i,\sum\limits_{i=0}^{n-1}\alpha^1_{i} p^i, \dots,
\sum\limits_{i=0}^{n-1}\alpha^{k-1}_{i} p^i).$$

A function $F:$ $\mathbb{Z}^{k}_{p}$ $\mapsto$ $\mathbb{Z}^{k}_{p}$
is said to satisfy the 1-lipschitz condition if $$|F(x) - F (y)|_p
\leq |x-y|_p, \text{ i.e. } F(x) \equiv F(x\text{ mod } p^{n})
\text{ mod } p^{n},$$  for every  $x, y \in \mathbb{Z}^{k}_{p}.$

Given a 1-lipschitz function $F:$ $\mathbb{Z}^{k}_{p}$ $\mapsto$
$\mathbb{Z}^{k}_{p}$, a mapping $F \text{ mod } p^n:\; x \mapsto
F(x) \text{ mod } p^n$, where $x \in (\mathbb{Z}/p^n\mathbb{Z})^k,$
is a well-defined mapping of the Cartesian power of the residue ring
modulo $p^n$ into itself. We call this mapping an induced function
modulo $p^n$.

A 1-lipschitz function $G:$ $\mathbb{Z}^{k}_{p}$ $\mapsto$
$\mathbb{Z}^{k}_{p}$ is transitive modulo $p^n$, if the map $G$ mod
$p^n$ is a permutation with a single cycle.
\section{Main result}

Everywhere we consider normalized Haar measure-preserving maps.

Take any 1-lipschitz measure-preserving map $F:\mathbb{Z}^{k}_{p}
\mapsto \mathbb{Z}^{k}_{p},$ where $k>1\in\mathbb{N},$ and represent
the entire set $\mathbb{Z}^{k}_{p}$ as a partition of subsets
$$\mathcal{F}_k(x_0) = \{x_0, F(x_0), \ldots, F^{p^k-1}(x_0)\},$$ where $x_0\in\mathbb Z^k_p:
x_0\equiv(0 ,\ldots, 0)\text{ mod } p.$ Subsets $\mathcal{F}_k(x_0)$
are disjoint for different $x_0$, as measure-preserving means
bijection [1]. And, as $F$ is 1-lipschitz, it's true that
$$F^{p^k}(x_{0}) = x,$$ $\text{ where } x\equiv (0,\;\dots,\;0) \text{ mod } p$, see definition of
1-lipschitz map in Section 2.

Take any 1-lipschitz measure-preserving map $G:\mathbb{Z}_{p}
\mapsto \mathbb{Z}_{p}$ and represent the entire set
$\mathbb{Z}_{p}$ as a partition of subsets
$$\mathcal{G}_{k}(y_0) = \{y_0, G(y_0), \ldots, G^{p^k-1}(y_0)\},$$
 where $y_0\in\mathbb Z_p: y_0\equiv 0 \text{ mod } p,$ $k$ is fixed natural. Proof, that it's partition of $\mathbb{Z}_{p}$,  is the same as for
$\mathcal{F}_{k}$.

Let $D :\mathbb{Z}^{m}_{p} \mapsto
\mathbb{Z}^{m}_{p},\;m\in\mathbb{N},$ be 1-lipschitz
measure-preserving map. $k\geq m$ is a fixed natural. Bijection
$T_{k,\;P}$, which can be associated with a permutation $P$ of $\{1,
\dots, p^k\}$, is defined for any $x \in\mathbb Z^m_p$, where $x \in
\mathcal{D}_m(x_0)$ and $D^j(x_0) = x$, $j=0,\dots,p^k-1,$ as
follows
$$T_{k,\;P}\circ D(x) = D^{P(j+1)}(x_0).$$

Define the map $H_{k}: \mathbb{Z}^k_p \mapsto \mathbb{Z}_p$ for any
natural $k>1$ as follows
$$H_{k}(\sum\limits_{i=0}^{\infty}\alpha^0_i
p^i,\sum\limits_{i=0}^{\infty}\alpha^1_{i} p^i, \dots,
\sum\limits_{i=0}^{\infty}\alpha^{k-1}_{i} p^i) =
\sum\limits_{i=0}^{\infty} \sum\limits_{j=0}^{k-1}\alpha^j_{i}
p^{ik+j},$$ where $\alpha^j_i\in\{0,\;\dots,\;p-1\}.$

\begin{theorem} For any 1-lipschitz measure-preserving transitive modulo $p$ map $F:\mathbb{Z}^{k}_{p} \mapsto
\mathbb{Z}^{k}_{p},$ where $k>1\in\mathbb{N},$ there are 1-lipschitz
measure-preserving transitive modulo $p^k$ map $G:\mathbb{Z}_{p}
\mapsto \mathbb{Z}_{p}$ and permutation $P$ of $\{1, \dots, p^k\}$,
that
$$G = H_{k} \circ T_{k,\;P}\circ F\circ H^{-1}_{k} \text{ and } F =
H^{-1}_{k} \circ T_{k,\;P^{-1}}\circ G\circ H_{k}.$$ Moreover, $F$
is ergodic iff $G$ is ergodic.\end{theorem}

\textit{Proof.} We can consider $(\sum\limits_{i=0}^{n-1}\alpha^0_i
p^i,\sum\limits_{i=0}^{n-1}\alpha^1_{i} p^i, \dots,
\sum\limits_{i=0}^{n-1}\alpha^{m-1}_{i} p^i) \in
(\mathbb{Z}/p^n\mathbb{Z})^m$ as element from $\mathbb{Z}^{m}_{p}$
for any $n, m \in \mathbb{N}$.

Take an arbitrary 1-lipschitz measure-preserving transitive modulo
$p^k$ function $G_1:\mathbb{Z}_{p} \mapsto \mathbb{Z}_{p}$. Take
such permutation $P$ of $\{1, \dots, p^k\}$, that
$$ H_{k}\circ T_{k,\;P}\circ F\circ H^{-1}_{k} \text{ mod } p^k = G_1 \text{ mod } p^k.$$

Consider $\hat{G}_n = H_{k}\circ T_{k,\;P}\circ F \circ H^{-1}_{k}
\text{ mod } p^{kn}$ as mapping on $\mathbb{Z}/p^{kn}\mathbb{Z}$ for
any natural $n$. There is a 1-lipschitz measure-preserving function
$G_n:\mathbb{Z}_{p} \mapsto \mathbb{Z}_{p}$, that $\hat{G}_n = G_n
\text{ mod } p^{kn}.$ It's true for $n=1$, see beginning of the
proof.

$n\rightarrow n+1$. $\hat{G}_{n+1}(x) \equiv \hat{G}_n(x\text{ mod }
p^{kn})\text{ mod } p^{kn}$ for any natural $n$, as $F$ is a
1-lipschitz map, see also definition of $H_k$ and $T_{k,\;P}$ in
this Section. We can describe $\hat{G}_n$ through a 1-lipschitz
measure-preserving function $G_n:\mathbb{Z}_{p} \mapsto
\mathbb{Z}_{p}$ by induction hypothesis. There is a 1-lipschitz
measure-preserving mapping $G_{n+1}:\mathbb{Z}_{p} \mapsto
\mathbb{Z}_{p}$, that
$$G_{n+1} \text{ mod } p^{kn} = G_n\text{ mod } p^{kn} \text{ and }
\hat{G}_{n+1} = G_{n+1} \text{ mod } p^{k(n+1)},$$ as mappings from
set of 1-lipschitz measure-preserving functions on $\mathbb{Z}_{p}$,
which equal $G_n$ modulo $p^{kn}$, take all possible distribution of
senior $k$ digits in the base $p$ system, if we consider them modulo
$p^{k(n+1)}$. $$G =\lim\limits_{n\rightarrow\infty} G_n.$$

The limit exists, as in the algebraic approach a p-adic integer is a
sequence $(a_n)_{n\geq1}$ such that $a_n$ is in
$\mathbb{Z}/p^n\mathbb{Z}$, and if $n \leq m$, then $a_n \equiv a_m
\text{ mod } p^n$.

We have already proofed that $G$ is a 1-lipschitz map, as $$G \text{
mod } p^{km+l} = G_{k(m+1)} \text{ mod } p^{km+l},\;m \in
\mathbb{N},\; 0\leq l \leq k-1,$$ $G_{k(m+1)}$ is a 1-lipschitz map.

$G$ is a measure-preserving map, as a 1-lipschitz function $D:
\mathbb{Z}^{k}_{p} \mapsto \mathbb{Z}^{k}_{p},$ is
measure-preserving if and only if it is bijective modulo $p^n$ for
any natural $n$, see Theorem 4.23 from [1]. Bijection modulo
$p^{k(n+1)}$ means bijection modulo $p^{m}$, where $m\leq k(n+1)$,
for any 1-lipschitz map, see definition of 1-lipschitz map in
Section 2. And $G_n$ is bijective modulo $p^{kn}$, as $\hat{G}_n =
H_{k}\circ T_{k,\;P}\circ F \circ H^{-1}_{k} \text{ mod } p^{kn}$ is
bijective.

Consider $F_n = H^{-1}_{k} \circ T^{-1}_{k,\;P}\circ G\circ H_{k}
\text{ mod } p^{n}$ as mapping on $(\mathbb{Z}/p^n\mathbb{Z})^k$ for
any $n \in \mathbb{N}$. We obtain 1-lipschitz measure-preserving map
$\tilde{F} = \lim\limits_{n\rightarrow\infty} F_n$ by the same
arguments.

Assumption $F\neq \tilde{F}$ implies the existence of $x \in
\mathbb{Z}^{k}_{p}$: $F(x)\neq \tilde{F}(x)$. It means, that there
is natural $m$: $F(x)\neq \tilde{F}(x) \text{ mod } p^m$. Hence, it
leads to the contradiction with $F\text{ mod } p^{n} = F_n$ for any
natural $n$, as $$F_n = H^{-1}_{k} \circ T^{-1}_{k,\;P}\circ G\circ
H_{k} \text{ mod } p^{n},$$
$$G \text{ mod } p^{kn} = H_{k}\circ T_{k,\;P}\circ F \circ H^{-1}_{k}
\text{ mod } p^{kn} = T_{k,\;P} \circ H_{k}\circ F \circ H^{-1}_{k}
\text{ mod } p^{kn}.$$ The last equitation is true, because it is
not important when we make the permutation: before or after
transformation $H_k$.

According to the preceding arguments, it can be shown that $F =
H^{-1}_{k} \circ T_{k,\;P^{-1}}\circ G\circ H_{k}$ and,
respectively, $G = H_{k} \circ T_{k,\;P}\circ F\circ H^{-1}_{k}$.
It's obvious, that $T^{-1}_{k,\;P} = T_{k,\;P^{-1}}$.

We are now to prove that transformations $T_{k,\;P}$ and $H_k$
preserve ergodicity.

A 1-lipschitz measure-preserving function $D: \mathbb{Z}^{k}_{p}
\mapsto \mathbb{Z}^{k}_{p}$ is ergodic if and only if $D$ is
transitive modulo $p^n$ for any natural $n$, see Proposition 4.35
from [1].

Necessity follows from the fact that if $F$ is transitive modulo
$p^n$ for every $n \in \mathbb{N}$, then $G$ is transitive modulo
$p^{kn}$.

Indeed, use $F = H^{-1}_{k} \circ T_{k,\;P^{-1}}\circ G\circ H_{k}$.
As an arbitrary permutation does not affect the property of
transitivity modulo $p^n$, $T_{k,\;P} \circ H_k \circ F \circ
H^{-1}_{k}\text{ mod } p^{kn}$ define the permutation with a single
cycle, as $F$ is ergodic. Therefore, $G \text{ mod } p^{kn}$ define
the same permutation.

Transitivity modulo $p^{k(n+1)}$ means transitivity modulo $p^{m}$,
where $m\leq k(n+1)$, for any 1-lipschitz map, see definition of
1-lipschitz map in Section 2.

Sufficiency follows from the fact, that we can repeat all our
arguments in the opposite direction, using the equality  $G = H_{k}
\circ T_{k,\;P}\circ F\circ H^{-1}_{k}$.$_\square$

\end{document}